\documentclass[12pt]{amsart}

\usepackage{amsmath}
\usepackage{amssymb}
\usepackage{amsthm}
\usepackage{graphicx}

\newtheorem{theorem}{Theorem}

\def\Z{\mathbb{Z}}
\def\eps{\varepsilon}

\begin{document}

\title{Computation of the Least Primitive root}

\author{Kevin J.~McGown}
\address{Department of Mathematics and Statistics,
  California State University at Chico}
\email{kmcgown@csuchico.edu}
\author{Jonathan P.~Sorenson}
\address{Computer Science and Software Engineering Department,
  Butler University, Indianapolis, IN 46208 USA}
\email{sorenson@butler.edu}

\date{\today}

\begin{abstract}
  Let $g(p)$ denote the least primitive root modulo $p$,
    and $h(p)$ the least primitive root modulo $p^2$.
  We computed $g(p)$ and $h(p)$ for all primes $p\le 10^{16}$.
  As a consequence we are able to prove that
    $g(p)<p^{5/8}$ for all primes $p>3$ 
    and that $h(p)<p^{2/3}$ for all primes $p$.
  We also present some additional results.
\end{abstract}

\maketitle

\section{Introduction and statement of results}\label{S:intro}

Let $g(p)$ denote the least primitive root modulo $p$, and
$h(p)$ denote the least primitive root modulo $p^2$.
We calculated $g(p)$ for all $p\leq 10^{16}$ and checked the condition
$g(p)=h(p)$.  This allows us to extend several known results.
Our computation took a total of roughly 645 days, wall time, 
on a linux cluster with 192 cores, or roughly 3 million CPU hours.

We computed $g(p)$ for all primes $p\leq 10^{16}$ and stored the values 
where $g(p)\geq 100$.  Here are the record values of $g(p)$:

\begin{tabular}{c|c}
$p$ & $g(p)$\\
\hline
$3$&$2$\\
$7$&$3$\\
$23$&$5$\\
$41$&$6$\\
$71$&$7$\\
$191$&$19$\\
$409$&$21$\\
$2161$&$23$\\
$5881$&$31$\\
$36721$&$37$\\
$55441$&$38$\\
$71761$&$44$\\
\end{tabular}\hspace{4ex}
\begin{tabular}{c|c}
$p$ & $g(p)$\\
\hline
$110881$&$69$\\
$760321$&$73$\\
$5109721$&$94$\\
$17551561$&$97$\\
$29418841$&$101$\\
$33358081$&$107$\\
$45024841$&$111$\\
$90441961$&$113$\\
$184254841$&$127$\\
$324013369$&$137$\\
$831143041$&$151$\\
$1685283601$&$164$\\
\end{tabular}\hspace{4ex}
\begin{tabular}{c|c}
$p$ & $g(p)$\\
\hline
$6064561441$&$179$\\
$7111268641$&$194$\\
$9470788801$&$197$\\
$28725635761$&$227$\\
$108709927561$&$229$\\
$386681163961$&$263$\\
$1990614824641$&$281$\\
$44384069747161$&$293$\\
$89637484042681$&$335$\\
$358973066123281$&$347$\\
$2069304073407481$&$359$\\
$4986561061454281$&$401$\\
$6525032504501281$&$417$\\
\end{tabular}

In order to show that the principal congruence subgroup $\Gamma(p)$ can be generated by
$[1, p; 0, p]$ and $p(p-1)(p+1)/12$ hyperbolic elements,
and construct such generators,
Grosswald employed the hypothesis that $g(p)<\sqrt{p}-2$ for all primes $p>409$, which is
now referred to as Grosswald's conjecture (see~\cite{MR49937}).
We remark that this conjecture is clearly true for $p$ large enough since, by
the Burgess inequality, one knows that $g(p)\ll p^{\frac{1}{4}+\eps}$;
however, explicit results where the upper bound is anywhere close to $p^{\frac{1}{4}}$ are difficult to come by,
except when $p$ is ridiculously large.
To illustrate this, we remark that in a subsequent paper, Grosswald showed that $g(p)\leq p^{0.499}$ when
$p>e^{e^{24}}\approx 10^{10^{10}}$ (see~\cite{MR636957}).

Nonetheless, in~\cite{MR3460815} it was shown that Grosswald's conjecture holds except possibly when
$p\in(2.5\cdot 10^{15}, 3.4\cdot 10^{71})$.  In~\cite{MR4055934} the bound of
$3.4\cdot 10^{71}$ was reduced to $10^{56}$,
and in~\cite{MR3584569} the conjecture was proved completely under the Generalized Riemann Hypothesis.
Together with these results, the list of records given in the table above allows one to easily verify that:

\begin{theorem}
Grosswald's conjecture holds except possibly when $p\in(10^{16},10^{56})$.
\end{theorem}

Previously it was verified in~\cite{MR2476579} that $g(p)=h(p)$ for all primes $p\leq 10^{12}$
except for $p=40487$ and $p=6692367337$.  In~\cite{MR4125902} it was pointed out that there are likely
infinitely many such primes $p$ for the simple reason that a lift of a generator of
$(\Z/p\Z)^\times$ to $(\Z/p^2\Z)^\times$ has a $(p-1)/p$ chance of being a generator.
However, we have verified that there are no additional examples where $g(p)\neq h(p)$ with
$p\leq 10^{16}$.

\begin{theorem}
The only primes $p\leq 10^{16}$ satisfying $g(p)\neq h(p)$ are
$40487$ and $6692367337$.
\end{theorem}

Perhaps this is not surprising.  If one believed the heuristic in the previous paragraph
and $x_n$ is the $n$-th odd prime $p$ where $g(p)\neq h(p)$, then one might expect the chances that
$g(p)=h(p)$ for all primes $p$ with $x_n<p\leq x$ to be roughly
$\prod_{x_n<p\leq x}(1-1/p)$;
denoting $x_1=40487$, $x_2=6692367337$,
this would lead one to expect that $x_3$ might lie
between $10^{19}$ and $10^{20}$.

As was the motivation for this work,
our computation complements the results in~\cite{ChenRoots,MR4055934,MR4125902},
leading to the following:

\begin{theorem}\label{T:3}
One has $g(p)<p^{5/8}$ for all primes $p>3$.
\end{theorem}

\begin{theorem}\label{T:4}
One has $h(p)<p^{2/3}$ for all primes $p$.
\end{theorem}

Previously, Pretorius \cite{pretorius2018smallest} showed that
  $g(p)<p^{0.68}$ for all primes $p$,
  and Chen \cite{ChenRoots} showed that $h(p)<p^{0.74}$ for all primes $p$.


%

In Section~2 we describe the methods used in our computations with a few
comments on their asymptotic time complexity.
In Section~3 we describe some auxillary programs used in the proofs
of Theorems~\ref{T:3} and~\ref{T:4}.
Finally, in Section~\ref{S:additional} we give some additional computational results.

\section{Computational methods} \label{S:computation}

We set a modulus $M=2\cdot3\cdots 29=6469693230$, and a core loop
  stepped through the integers $a\bmod M$ with $\gcd(a,M)=1$.
The $\phi(M)=1021870080$ values of $a$ were shared evenly among
  192 processor cores and processed independently.

For each such $a$, the primes in the arithmetic progression
  $a\bmod M$ up to $10^{16}$ are found using the sieve of Eratosthenes.
This requires sieving an interval of size only $10^{16}/M\approx 1.5\times10^6$,
  which fits in cache.
Since $\pi(10^8)-10$, the number of primes used to sieve the interval, is
  5761445, of a similar magnitude, the work is reasonably balanced.

In addition to this, the progression $a-1 \bmod M$ is sieved to generate
  complete factorizations, so that for each prime $p\equiv a\bmod M$ we
  find, we have handy the complete factorization of $p-1$.
With this information, we then find $g(p)$ and check if $g(p)=g(p^2)$ 
  using Montgomery multiplication \cite{MR777282}
  to reduce the cost of modular exponentiations.
For each of the integers $n\equiv a-1 \bmod M$ that we factor completely,
  the odd prime divisors of $n$ up to $10^8$ are packed, in binary, into
  a single 64-bit integer to save space.

The 192-core cluster at Butler University runs linux;
  the code was written in C++ with MPI
  to manage parallel execution and communication.
The total wall time for the computation was 645 days, with two
  brief interruptions requiring restarts from saved checkpoints.

If we generalize our method to an algorithm to find $g(p)$ and $h(p)$
  for all primes $p\le n$,
  our algorithm takes $O(n)$ arithmetic operations under the assumption that
  the average values of $g(p)$ and $h(p)$ are bounded.
(See, for example, \cite{EM97}.)
We assume that $M$, our product of small primes,
  satisfies $n^{1/3}\le M \le n^{2/3}$, say, so that
  $\phi(M)/M \approx \log\log n$.
This serves the same role as a wheel in prime sieves, so that
  finding all primes up to $n$ and factoring all integers
  in the $a-1 \bmod M$ residue classes with $\gcd(a,M)=1$
  take $O(n)$ operations each.
(See \cite{Sorenson2015} for references on prime sieves.)
Then, for each prime, we compute $g(p)$ and $h(p)$ which takes
  $O( h(p) \log p )$ arithmetic operations for each $p$,
  as a modular exponentiation takes $O(\log q)$ operations
  if the exponent is bounded by $q$.
By the prime number theorem and our heuristic assumption that the
  average value of $h(p)$ is bounded, we get $O(n)$ operations here as well.

%
%

\section{Proof of Theorems~\ref{T:3} and~\ref{T:4}}\label{S:proofs}

Our proofs for both theorems are computational, and used a total of five
  programs.
At a high level, our approach is not dissimilar to what Chen outlines
  in Remark 2.6 in \cite{ChenRoots}.

\subsection{Theorem \ref{T:3}}
To prove Theorem~\ref{T:3} 
we used Theorem~3 from \cite{MR4055934}, namely that if
$$
H^2 - \frac{\pi^2}{6} \frac{B(H/h)^{2r-1}}{A(H/h)^{2r}}
  \left\{
	  \left( 2+\frac{s-1}{\delta}\right) 2^{\omega-s}
	  \right\}^{2r}
  h \sqrt{p} W(p,h,r)>0
$$
then $g(p)<H$.
Following the notation from that paper,
  we set $\alpha=5/8$, $H=p^\alpha$, $h=\lceil 2.01\cdot p^{2\alpha-1}\rceil$,
  and $\omega:=\omega(p-1)$, the number of prime divisors of $p-1$.
If we set $e$ to be an even divisor of $p-1$, 
  then $s=\omega(p-1)-\omega(e)$, the number
  of distinct primes dividing $p-1$ that do not divide $e$.
Let $p_1, p_2, \ldots, p_s$ be theses primes.
Then $\delta:=1-\sum_{i=1}^s 1/p_i$.
We also have
\begin{eqnarray*}
	W(p,h,r) &\le&
	\begin{cases}
	\sqrt{2p} \left( \frac{2r}{eh} \right)^r+2r-1
	&
	\mbox{for $r\ge 1$\,,} \\
	3\left(1+\frac{\sqrt{p}}{h^2}\right)
	&
	\mbox{for $r=2$\,.}
	\end{cases}
\end{eqnarray*}

Our first program, written in C++, works to narrow the range of $p$
  for when we cannot yet prove that $g(p)<p^\alpha$.
We begin by knowing from our computation
  above and the previous result that $g(p)<p^{5/8}$ for $p>10^{22}$
  (Corollary 1 from \cite{MR4055934})
  that we need only consider primes $p$ with $10^{16} < p < 10^{22}$.
We then tried a range of $r$ values from $1$ to $10$,
  and bounded $\omega$ explicitly (Theorem~11 from \cite{omegabound}),
  and set $s=0$ (which gives $\delta=1$) to lower the upper bound on $p$,
  using bisection.
Next, with a somewhat narrower range on $p$, this gives us a range for
  $\omega$, from $2$ to $18$ in this case.
We handled each value of $\omega$ as a separate case.
For each $\omega$, we can compute an optimal value of $s$ to use,
  by minimizing $(2+(s-1)/\delta)2^{\omega-s}$, as the value of
  $\omega-s$ can be used to give a lower bound on $\delta$.
With these more refined parameters, we closed the gap on $p$ values
  except for when $10\le \omega \le 15$.

To deal with the remaining values of $\omega$ we employed a tree algorithm,
similar what what is described in Section~3.1 of~\cite{MR3584569}.
Our version of the tree algorithm is to set a mask with a fixed number of
  bits -- 14 in this particular case.
A bit of the mask was $1$ if the prime corresponding to that bit position
  divided $p-1$.
For example, the mask 1010011 would mean the primes $2,3,11,17$ divide $p-1$
  and no other primes up to $43$ (the 14th prime) divide $p-1$.
This extra information allows a better lower bound on $\delta$, and in
  all cases we found that setting $s=\omega-1$ worked.
Looping through all odd-valued masks up to $2^{14}$,
  we narrowed the gaps on $p$ significantly, and in cases where a gap
  remained, we wrote the upper and lower limits, and the mask information
  to a text file to be read by a second program.
If the mask value had enough bits set so that all primes dividing $p-1$ were
  known, the information was written to a different text file for our
  third program to process.

This first program takes under a minute to run.

The second program read in upper and lower bounds, and a mask to use
  to create a modulus to create an arithmetic progression to sieve
  for primes.
For each prime found, its primitive root was found and checked.
Processing all the work generated by the first program took about five hours.
This second program was a SageMath script.
(If it had taken longer, we would have rewritten it in C++.)

The third program was given a complete list of prime divisors of $p-1$,
  so it enumerated all possible exponents on each prime specified by
  the mask to create candidate primes $p$ that were prime tested
  (see Theorem~4.1.1 from \cite{CrandallPomerance})
  and then primitive roots were found and checked.
This program took only a few hours to do the work it was given by
  the first program.

We reran the first program multiple times, trying different mask sizes
  and estimating the amount of work it produced for the other
  two programs, and we found 14 to be the best choice.

\subsection{Theorem \ref{T:4}}
Our proof of Theorem~\ref{T:4} 
is based on the recent work of Chen \cite{ChenRoots},
and consists of two programs that work in roughly the same way as above.

The first program is based on 
equation~(2.5) from \cite{ChenRoots}:
$$
  \delta \frac{\phi(e)}{e} \left\{
	  p^\alpha - 
	  \left(2+ \frac{s-1}{\delta}\right) 2^{\omega(e)}
	  A(p) \sqrt{p}\log p
	  \right\}
  -\sqrt{5} \pi p^{\alpha-1/4} > 0.
$$
Note that although the $A(p)$ function here is from \cite{ChenRoots} and is not
  the same as the $A(H/h)$ function in the previous section,
  the other parameters $e, s, \delta$ are defined similarly.
Our proof used $\alpha=0.661<2/3$, so what we proved is slightly
  stronger than the statement of the theorem.
(We found that using $0.660$ would have required more work
  than is feasible -- $0.661$ seems the limit of this method for $h(p)$.)
We required an explicit lower bound for $\phi(e)/e$ which can be found
  in \cite{RosserSchoenfeld}.

Otherwise, the overall program structure is basically the same as for
  the first program for the proof of Theorem~\ref{T:3}.

We had no upper bound on $p$ for our proof, so computing one was the first step,
  using explicit bounds for $\omega$ and $\phi(e)/e$, giving $10^{233}$.
This meant $\omega\le 105$.
Proving $\omega=2$ to $8$ and $17$ to $105$ was straightforward, but for
  $\omega$ between $9$ and $16$ inclusive we used a mask
  of 16 bits and handled cases as before.
Note that $s=0$ worked for $\omega$ up to $6$, then $s=4,6$ for
  $\omega=7,8$ respectively, and we used $s=14$ for $\omega\ge17$.
We used $s=\omega-2$ in the troublesome range $9$--$16$.

The only change in the second program was to compute $h(p)$ instead of $g(p)$.
Like the previous subsection, it was a SageMath script as well.
We did not use a third program as in the previous subsection -- the second program was able to do the work instead.

\section{Additional computational results}\label{S:additional}

\begin{figure}
\includegraphics[width=5in]{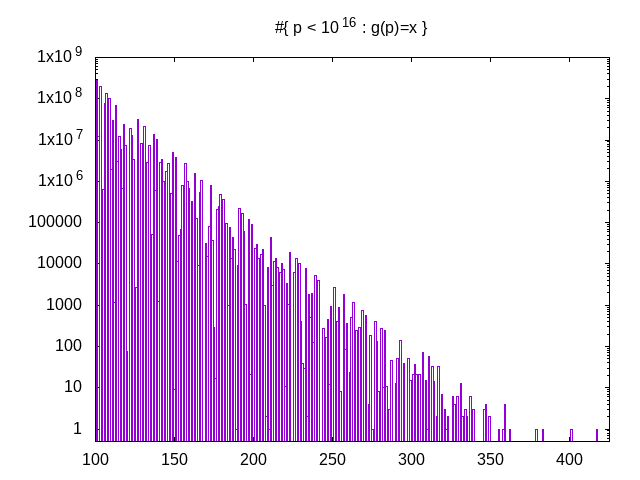}
\caption{Log-histogram of $g(p)$}
\label{F:hist}
\end{figure}

We felt it would be of interest to include a log-histogram of the values of $g(p)$,
where $g(p)>100$.  See Figure~\ref{F:hist}.

\begin{table}
\begin{tabular}{|r|c|c|c|}
	\hline
     $p$ & $g(p)$ & $F(p)$ & $g(p)/F(p)$\\
     \hline
          29418841 &	101 &	247.87 	&	0.407471\\
          33358081 &	107 	&250.961 &	0.426361\\
          45024841 &	111 	&258.388 &	0.429586\\
          90441961 &	113 	&275.932 &	0.409521\\
         184254841 &	127 	&294.212 &	0.431661\\
         324013369 &	137 	&308.979 &	0.443396\\
         831143041 &	151 	&334.133 &	0.451916\\
        1685283601 &	164 	&353.416 &	0.464042\\
        6064561441 &	179 	&389.207 &	0.459909\\
        7111268641 &	194 	&393.733 &	0.492719\\
        9470788801 &	197 	&401.919 &	0.490148\\
       28725635761 &	227 	&434.111 &	0.522908\\
      108709927561 &	229 	&473.726 &	0.483402\\
      386681163961 &	263 	&512.476 &	0.513195\\
     1990614824641 &	281 &	563.874& 	0.498338\\
    44384069747161 &	293 	&665.223& 	0.440453\\
    89637484042681 &	335 &	688.859 &	0.486311\\
   358973066123281& 	347 &	736.23& 		0.47132\\
  2069304073407481 &	359 &	797.351& 	0.450241\\
  4986561061454281 &	401 &	828.577 &	0.483962\\
  6525032504501281 &	417 &	838.194 &	0.497498\\ \hline
\end{tabular}
\caption{Growth of record values of $g(p)$}
\label{T:gF}
\end{table}

Let $\hat{g}(p)$ denote the least prime primitive root.
(Clearly $g(p)\leq \hat{g}(p)$, but they may not be equal.)
In~\cite{MR1433261}, based on a probabilistic model, it was conjectured that
$$
\limsup_{p\to\infty}\frac{\hat{g}(p)}{e^\gamma \log p(\log \log p)^2}=1
$$
and numerical evidence to support this was given using a table of record values for 
$\hat{g}(p)$, $p<2^{31}\approx 2.1\cdot 10^9$.  We carry out a similar computation
for $g(p)$ using the record values given in Section~\ref{S:intro}.
Set $F(p)=e^\gamma \log p(\log \log p)^2$.
We certainly do not expect $\limsup g(p)/F(p)$ to equal one, but just for sake
of comparison we compute $g(p)/F(p)$ as p ranges over the record values.
See Table~\ref{T:gF}.

The value of the ratio is somewhat close to $1/2$, but we dare not speculate 
soley based on this numerical data.

\begin{figure}[hb]
\includegraphics[width=5in]{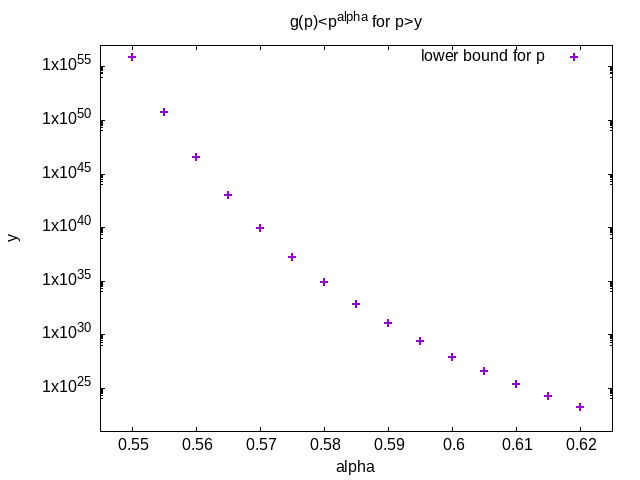}
\caption{Graph representation of Theorem~\ref{T:alphas}}
\label{F:alpha}
\end{figure}

We conclude with the following theorem,
which is a refinement of Theorem~\ref{T:3},
proven using
methods similar to those discussed in Section~\ref{S:proofs}.

\begin{theorem}\label{T:alphas}
	If $p>p_\alpha$ then $g(p)<p^\alpha$ for the values of $\alpha$ and
	$p_\alpha$ in the following table:
	\begin{quote}
\begin{tabular}{|l|c||l|c|} \hline
	$\alpha$ & $p_\alpha$ &
	$\alpha$ & $p_\alpha$ \\ \hline
	$0.62$ & $1.91\times10^{23}$ &
	  $0.58$ & $8.21\times10^{34}$ \\
	$0.615$ & $1.97\times10^{24}$&
	  $0.575$ & $1.77\times10^{37}$ \\
	$0.61$ & $2.47\times10^{25}$ &
	  $0.57$ & $7.98\times10^{39}$ \\
	$0.605$ & $4.00\times10^{26}$&
	  $0.565$ & $9.42\times10^{42}$ \\
	$0.6 $ & $8.52\times10^{27}$&
	  $0.56$ & $3.57\times10^{46}$ \\
	$0.595$ & $2.53\times10^{29}$&
	  $0.555$ & $6.08\times10^{50}$ \\
	$0.59$ & $1.09\times10^{31}$&
	  $0.55$ & $7.27\times10^{55}$ \\
	$0.585$ & $7.24\times10^{32}$ && \\
\hline
\end{tabular}
	\end{quote}
\end{theorem}
The graph in Figure~\ref{F:alpha} illustrates the table above.

\nocite{*}
\bibliographystyle{plain}
\bibliography{primitive_root}

\end{document}